\documentclass[11pt]{article}
\usepackage{}
\usepackage[total={6in, 9in}]{geometry}
 \usepackage{amsfonts}
\usepackage{amsthm}
\usepackage{amssymb}
\usepackage{amsmath,amsfonts}
\usepackage{mathrsfs}
\usepackage{cases}
\usepackage{latexsym,bm}
\usepackage{indentfirst}
\usepackage{color}
\usepackage{ifpdf}
\usepackage{graphicx}
\usepackage{psfrag}
\usepackage{dsfont}
\usepackage[
pdfauthor={CESS},
pdftitle={Spanning trails},
pdfstartview=XYZ,
bookmarks=true,
colorlinks=true,
linkcolor=blue,
urlcolor=blue,
citecolor=blue,
bookmarks=true,
linktocpage=true,
hyperindex=true
]{hyperref}

\newtheorem{THM}{\textbf{Theorem}}

\newtheorem{LEM}{\textbf{Lemma}}

\newcommand{\pf}{\noindent\textbf{Proof}.\quad}

\newcommand{\CC}{\mathcal{C}}

\newcommand{\pbar}{\overline{\varphi}}
\newcommand{\psbar}{\overline{\psi}}

\begin{document}
\title{A note on Gupta's co-density conjecture}
\author{
	Guantao  Chen\\
	Georgia State University, Atlanta, GA 30302\\
	{\tt gcehn@gsu.edu}
	\and
Songling Shan\\
	Illinois State  University, Normal, IL 61790\\
		{\tt sshan12@ilstu.edu}
}

\date{\today}
\maketitle

\emph{\textbf{Abstract}.}
Let $G$ be a multigraph. A subset $F$ of $E(G)$ is an edge cover of $G$ if every vertex of $G$
is incident to an edge of $F$. The cover index, $\xi(G)$, is the largest number of edge covers into which 
the edges of $G$ can be partitioned.  Clearly $\xi(G) \le \delta(G)$, the minimum degree of $G$. 
For $U\subseteq V(G)$, denote by $E^+(U)$ the set of edges incident to
a vertex of $U$. When $|U|$ is odd, to cover all the vertices of $U$, any edge cover needs to contain at least $(|U|+1)/2$ edges from $E^+(U)$, indicating 
$ \xi(G) \le |E^+(U)|/ (|U|+1)/2$. 
Let $\rho_c(G)$, the co-density of $G$, be defined as  the minimum of $|E^+(U)|/((|U|+1)/2)$ ranging over all $U\subseteq V(G)$ with $|U| $  odd 
and at least 3.   Then  $\rho_c(G)$ provides another upper bound on $\xi(G)$.  Thus $\xi(G) \le \min\{\delta(G), \lfloor \rho_c(G) \rfloor \}$. 
For a lower bound on $\xi(G)$, 
in 1967, Gupta conjectured that $\xi(G)  \ge  \min\{\delta(G)-1, \lfloor \rho_c(G) \rfloor \}$. Gupta showed that the 
conjecture is true when $G$ is simple, and Cao et al. verified this conjecture when $\rho_c(G)$ is not an integer. 
In this note, we confirm the conjecture when the maximum multiplicity of $G$ is at most two or $ \min\{\delta(G)-1, \lfloor \rho_c(G) \rfloor \} \le 6$.

\emph{\textbf{Keywords}.} Edge cover; cover index; co-density; chromatic index.  

\vspace{2mm}

\section{Introduction}

Graphs in this paper may contain multiple edges but contain no loops.  
Let $G$ be a graph.
Denote by $V(G)$ and  $E(G)$ the vertex set and the  edge set of $G$,
respectively. For $v\in V(G)$, $d_G(v)$, the degree of $v$, is the number of 
edges of $G$ that are incident with $v$. For $S\subseteq V(G)$, 
the subgraph of $G$ induced on  $S$ is denoted by $G[S]$, and on $V(G)\setminus S$ is denoted by
$G-S$. For notational simplicity we write $G-x$ for $G-\{x\}$. For  $e\in E(G)$, $G-e$
is obtained from $G$ by deleting the edge $e$.  For an edge $e\not\in E(G)$, $G+e$ 
is obtained by adding the edge $e$ to $G$. 
Let $A,B \subseteq V(G)$ be disjoint.  We denote by 
$E_G(A)$ the set of edges with both endvertices in $A$,
$E_G(A,B)$ the set of edges with one endvertex  in $A$
and the other endvertex in $B$, and by $E_G^+(A)$ the set of edges of $G$
 incident with a vertex of $A$. Note that $E_G^+(A)$ is 
the union of $E_G(A)$ and $E_G(A,V(G)\setminus A)$. 
When $A=\{x\}$, we simply write $E_G(\{x\}, B)$ as $E_G(x,B)$. 
Let $e_G(A)=|E_G(A)|$, $e_G(A,B)=|E_G(A,B)|$, and $e_G^+(A)=|E_G^+(A)|$. 
When $G$ is clear from the context, we skip the subscript $G$
from the corresponding notation.

Let $F\subseteq E(G)$.  The set $F$  \emph{saturates} $v\in V(G)$ 
if $v$ is incident in $G$ with an edge from $F$; otherwise  $F$ \emph{misses} $v$. 
For $S\subseteq V(G)$, we say $F$ saturates $S$ if $F$ 
saturates every vertex of $S$. 
We call  $F$  \emph{an edge cover} of $G$   if $F$ saturates $V(G)$. The cover index, $\xi(G)$, is the largest number of edge covers into which 
the edges of $G$ can be partitioned.  Clearly $\xi(G) \le \delta(G)$, the minimum degree of $G$. 
For $U\subseteq V(G)$,  $|U|$ is odd,  to cover all the vertices of $U$, any edge cover needs to contain at least $(|U|+1)/2$ edges from $E^+(U)$.
Therefore, we have $ \xi(G) \le e^+(U)/ (|U|+1)/2$. 
Let $\rho_c(G)$, the \emph{co-density} of $G$, be defined as  the minimum of $e^+(U)/((|U|+1)/2)$ ranging over all $U\subseteq V(G)$ with $|U| $  odd 
and at least 3.   Then  $\rho_c(G)$ provides another upper bound on $\xi(G)$.  Thus $\xi(G) \le \min\{\delta(G), \lfloor \rho_c(G) \rfloor \}$. 
For a lower bound on $\xi(G)$, 
in 1967, Gupta~\cite{MR0491300} conjectured that $\xi(G)  \ge  \min\{\delta(G)-1, \lfloor \rho_c(G) \rfloor \}$,  and he proved the conjecture when $G$
is simple~\cite{Gupta}. 
 In 2019, Cao, the first author, Ding, Jing and Zang~\cite{1906.06458} verified this conjecture when $\rho_c(G)$ is not an integer. 
In this note, we generalize Gupta's result from simple graphs to graphs with maximum multiplicity at most two  and confirm the conjecture for  graphs $G$ with small $\delta(G)$
and $\rho_c(G)$ as stated below. 

\begin{THM}\label{thm:main}
	Let $G$ be a graph and $k=\min\{\delta(G)-1, \lfloor \rho_c(G) \rfloor \}$. If the maximum multiplicity of $G$
	is at most 2 or $k\le 6$, then 
 $G$ has at least $k$ edge-disjoint  edge covers. 
\end{THM}

The remainder of this paper is organized as follows. In the next section, we provide some notation and 
preliminaries; and in the last section, we prove Theorem~\ref{thm:main}. 

\section{Notation and Lemmas}

For two integers $p$ and $q$, let $[p,q]=\{i\in \mathbb{Z}:  p\le i\le q\}$. 
Let $G$ be a graph and $m\ge 0$ be an integer.  An  \emph{edge $m$-coloring} of $G$ is a map $\varphi$:  $E(G) \rightarrow [1,m]$   that assigns to every edge $e$ of $G$ a color $\varphi(e) \in [1,m]$  such that  no two adjacent edges receive the same color.  Denote by $\CC^m(G)$ the set of all edge $m$-colorings of $G$. 
The {\it chromatic index\/} $\chi'(G)$ is  the least integer $m\ge 0$ such that $\CC^m(G) \ne \emptyset$. An edge $e$ of $G$ is {\it critical} if 
$\chi'(G -e) < \chi'(G)$.

For a vertex  $v\in V(G)$ and  a coloring $\varphi \in \CC^{m}(G)$ for some integer $m\ge 0$, define the two color sets
$
\varphi(v)=\{\varphi(f)\,:\, \text{$f$ is incident to $v$ in $G$}\}$ and $ \pbar(v)=[1, m] \setminus\varphi(v).
$
We call $\varphi(v)$ the set of colors \emph{present} at $v$ and $\pbar(v)$
the set of colors \emph{missing} at $v$. 
For a color $\alpha$, the edge set $E_{\alpha} = \{ f\in E(G)\, |\, \varphi(f) = \alpha\}$ is called a {\it color class}. Clearly, 
$E_{\alpha}$ is a {\it matching} of $G$ (possibly empty). 
For two distinct colors $\alpha,\beta$,  the  subgraph of $G$
induced by $E_{\alpha}\cup E_{\beta}$ is a union of disjoint 
paths and  even cycles, which are referred to as   \emph{$(\alpha,\beta)$-chains} of $G$
with respect to $\varphi$.  
 For $x,y\in V(G)$, if $x$ and $y$
are contained in the same  $(\alpha,\beta)$-chain with respect to $\varphi$, we say $x$ 
and $y$ are \emph{$(\alpha,\beta)$-linked}.
Otherwise, they are \emph{$(\alpha,\beta)$-unlinked}. 

For a vertex $v$, let $C_v(\alpha, \beta, \varphi)$ denote the unique $(\alpha, \beta)$-chain 
containing $v$.  
If $C_v(\alpha, \beta, \varphi)$ is a path, we just write it as $P_v(\alpha, \beta, \varphi)$.
The notation $P_v(\alpha, \beta, \varphi)$ is commonly used when we know  $|\pbar(v)\cap \{\alpha,\beta\}|=1$.  If we interchange the colors $\alpha$ and $\beta$
on an $(\alpha,\beta)$-chain $C$ of $G$, we briefly say that the new coloring is obtained from $\varphi$ by an 
{\it $(\alpha,\beta)$-swap} on $C$, and we write it as  $\varphi/C$. 
This operation is called a \emph{Kempe-change}.

For an integer $s\ge 1$, 
a graph $G$ is   {\em $s$-dense }  if $|V(G)| \ge 3$ is odd and $|E(G)|=(|V(G)|-1)s/2$.  As a largest  matching in $G$
can have size at most $|V(G)|-1$, the lemma below is a consequence  of being $s$-dense, where a matching is \emph{near perfect} 
in $G$ if it misses only one vertex of $G$.

\begin{LEM}\label{k-dense subgraph}
	Let $G$ be an $s$-dense graph  with  $\chi'(G) = s$ for some integer $s\ge 1$, and let $\varphi\in \CC^{s}(G)$. 
Then for any  two distinct $u,v\in V(G)$,  we have $\pbar(u)\cap \pbar(v)=\emptyset$. In particular, 
	each color class of $\varphi$
	is a near perfect matching of $G$, and each vertex $v\in V(G)$ is missed by exactly  $s-d(v)$ of the color classes of $\varphi$. 
\end{LEM}

Let $\rho(G)$, the \emph{density} of $G$, be defined as  the  maximum of $e(U)/((|U|+1)/2)$ ranging over all $U\subseteq V(G)$ with $|U| $  odd 
and at least 3. 
In the 1970s,  Goldberg~\cite{MR354429} and Seymour~\cite{MR532981}
independently conjectured that  every graph $G$ satisfies $\chi'(G) \le \max\{\Delta(G)+1, \lceil \rho(G) \rceil\}$. 
Over the past four decades this conjecture has been a subject of extensive research.  In 2019, the first author, Jing, and Zang~\cite{1901.10316} announced a complete
proof of the Conjecture, and the lemma below is equivalent to the conjecture by Goldberg and Seymour.

\begin{LEM}\label{G-S-C}
	Let $G$ be a graph and $e\in E(G)$.  If $\chi'(G)=s+1\ge \Delta(G)+2$ and $\chi'(G-e)=s$, 
	then
	$G-e$ has an $s$-dense subgraph $H$ containing the endvertices of $e$ such that $e$ is also a  critical edge of $H+e$.
\end{LEM}

\begin{LEM}\label{lemma:special-coloring}
Let $G$ be a graph  and $k\ge 1$ be an integer. Suppose  $\Delta(G) \le k+1$ and  $\chi'(G)\le k+2$. 
Let $S$ be the set of vertices of $G$ with degree at most $k/2$. Then there exists $\varphi\in \CC^{k+2}(G)$
satisfying the following properties:
\begin{enumerate}
	\item  the color $k+2$ is missing at every vertex in $S$; 
	\item If $k+1\in \varphi(x)$ for some $x\in S$, then $P_x(k+1,k+2, \varphi)$ ends 
	at a vertex in $V(G)\setminus S$. 
\end{enumerate}
\end{LEM}

\pf For any  $\varphi\in \CC^{k+2}(G)$ we define 
\begin{eqnarray*}
s_\varphi &=&|\{x\in S: k+2\in \varphi(x)\}|, \quad \text{and} \\ 
c_\varphi&=&|\{P_x(k+1, k+2, \varphi):  \text{$P_x(k+1, k+2,\varphi)=P_y(k+1,k+2,\varphi)$ for  distinct $x, y\in S$}\}|, 
\end{eqnarray*}
to be respectively the number of vertices in $S$ that present  the color $k+2$
and the number of    $(k+1, k+2)$-chains (path-chain) with both endvertices in $S$ under $\varphi$. 
We choose $\varphi\in \CC^{k+2}(G)$ with $s_\varphi$ minimum and subject to this, with $c_\varphi$ minimum. 
If $s_\varphi=c_\varphi=0$, then we are done. Thus we assume  $s_\varphi+c_\varphi>0$.  

Consider first that $s_\varphi>0$. Let $x\in S$ such that $k+2\in \varphi(x)$. 
Since $d(x) \le k/2$, there exists $\alpha\in [1,k]$ such that $\alpha \in \pbar(x)$. 
We consider $P_x(\alpha,k+2,\varphi)$.  If $P_x(\alpha,k+2,\varphi)$ ends at a vertex not in $S$ or ends at a vertex from $S$
that presents $k+2$,  then $\psi:= \varphi/P_x(\alpha,k+2,\varphi)$  is an edge $(k+2)$-coloring of $G$ with $s_\psi <s_\varphi$. 
Thus we assume that $P_x(\alpha,k+2,\varphi)$ ends at a vertex $y\in S\setminus\{x\}$ such that $\alpha\in \varphi(y)$
and $k+2\in \pbar(y)$. Let 
$$
P_x(\alpha,k+2,\varphi)=v_0v_1\ldots v_{2t-1} v_{2t}, 
$$
for some integer $t\ge 1$, where $v_0:=x$ and $v_{2t}:=y$.  
 
Since $|\varphi(x)\cup \varphi(y)|\le d(x)+d(y) \le k$,   we have $ \pbar(x)\cap \pbar(y) =[1,k+2]\setminus  (\varphi(x)\cup \varphi(y))\ne \emptyset$.
Let $i\in [1,2t]$ be the smallest index such that $\pbar(v_i)\cap \pbar(x) \ne \emptyset$.  
As $k+2 \in \varphi(x)$, $k+2\not\in \pbar(v_i)\cap \pbar(x)$.  Among all the edge $(k+2)$-colorings  $\xi$
with $s_\xi=s_\varphi$, $c_\xi=c_\varphi$,   and $P_x(\alpha,k+2,\xi)=P_x(\alpha,k+2,\varphi)$, 
we may assume $\varphi$ is the one such that the index $i$ is smallest.

If $i=1$, then  simply recoloring  $xv_1$ by a color from $\pbar(v_1)\cap \pbar(x)$ gives a new coloring $\psi$
with $s_\psi <s_\varphi$.  Thus $i\ge 2$. 
Let $\beta \in \pbar(v_i)\cap \pbar(x) \subseteq [1,k+1]$. By the minimality of $i$, we have $\beta \in \varphi(v_{i-1})$. 
As $d(v_{i-1}) \le k+1$ and $\alpha,\beta, k+2\in \varphi(v_{i-1})$,  there exists $\gamma\in \pbar(v_{i-1}) \subseteq [1,k+2]\setminus \{\alpha, \beta, k+2\}$. 

If $v_i$ and $v_{i-1}$ are not $(\beta,\gamma)$-linked, then let $\psi$ be obtained 
by doing a  Kemple-change on $P_{v_i}(\beta, \gamma, \varphi)$ and then recoloring the edge $v_{i-1}v_i$ on $P_x(\alpha,k+2, \varphi)$ by $\gamma$. 
Note that $s_\psi \le s_\varphi$, 
and we have that $\alpha \in \psbar(v_{i-1})$ or $k+2 \in \psbar(v_{i-1})$, and $P_x(\alpha,k+2,\psi)=P_{v_{i-1}}(\alpha,k+2,\psi)$. 
If  $\alpha \in \psbar(v_{i-1})$, then we  can do a Kempe-change on  $P_x(\alpha,k+2,\psi)$ to decrease $s_\psi$ and so to decrease $s_\varphi$. 
Thus we assume that $k+2 \in \psbar(v_{i-1})$.  If $v_{i-1} \in S$, then we have $s_\psi <s_\varphi$ already. Thus 
we assume $v_{i-1} \notin S$.  Then we  can do a Kempe-change on  $P_x(\alpha,k+2,\psi)$ to decrease $s_\psi$ and so to decrease $s_\varphi$. 

 Thus we  assume  now that  $v_i$ and $v_{i-1}$ are  $(\beta,\gamma)$-linked.  Then let $\psi=\varphi/P_{v_i}(\beta,\gamma,\varphi)$. 
We have $s_\psi=s_\varphi$,  $c_\psi=c_\varphi$, and $P_x(\alpha,k+2,\psi)=P_y(\alpha,k+2,\varphi)$. 
However, we have $ \beta \in \psbar(v_{i-1})\cap \psbar(x) \ne \emptyset$, contradicting the choice of $\varphi$. 
Thus $s_\varphi=0$. 

Consider now that $s_\varphi=0$ but $c_\varphi>0$.  Then there exist distinct $x,y\in S$
such that $P_x(k+1,k+2,\varphi)=P_y(k+1,k+2,\varphi)$. Note that $k+2\in \pbar(x)\cap \pbar(y)$, 
and $P_x(k+1,k+2,\varphi)$ is internally disjoint from $S$ as  $s_\varphi=0$. 
Let 
$$
P_x(k+1,k+2,\varphi)=v_0v_1\ldots v_{2t} v_{2t+1}, 
$$
for some integer $t\ge 0$, where $v_0:=x$ and $v_{2t+1}:=y$.  
Since $d(x)+d(y) \le k$ and $k+1\in \varphi(x)\cap \varphi(y)$,   we have $ (\pbar(x)\cap \pbar(y))\cap [1,k] \ne \emptyset$. 

Let $i\in [1,2t+1]$ be the smallest index such that $(\pbar(v_i)\cap \pbar(x))\cap [1,k]  \ne \emptyset$.  
 Among all the edge $(k+2)$-colorings  $\xi$
with $s_\xi=0$, $c_\xi=c_\varphi$ and $P_x(k+1,k+2,\xi)=P_x(k+1,k+2,\varphi)$, 
we may assume $\varphi$ is the one such that the index $i$ is smallest.

If $i=1$, then we simply recolor $xv_1$ by a color from $(\pbar(v_1)\cap \pbar(x))\cap [1,k]$ gives a new coloring $\psi$
with $c_\psi <c_\varphi$.   Since the new color is from $[1,k]$, we still have  $s_\psi=0$. 
Thus $i\ge 2$. 
Let $\beta \in( \pbar(v_i)\cap \pbar(x))\cap [1,k]$. By the minimality of $i$, we have $\beta \in \varphi(v_{i-1})$. 
As $d(v_{i-1}) \le k+1$ and $\beta,  k+1, k+2\in \varphi(v_{i-1})$,  there exists $\gamma\in \pbar(v_{i-1}) \subseteq [1,k]\setminus \{\beta\}$. 

If $v_i$ and $v_{i-1}$ are not $(\beta,\gamma)$-linked, then let $\psi$ be obtained 
by doing a  Kemple-change on $P_{v_i}(\beta, \gamma, \varphi)$ and then recoloring the edge $v_{i-1}v_i$ on $P_x(k+1,k+2, \varphi)$ by $\gamma$. 
Then $c_\psi <c_\varphi$.     Since $\beta, \gamma\in [1,k]$, we still have $s_\varphi=0$.  Thus   $v_i$ and $v_{i-1}$ are  $(\beta,\gamma)$-linked.
Then let $\psi=\varphi/P_{v_i}(\beta,\gamma,\varphi)$. 
We have $s_\varphi=0$,  $c_\psi=c_\varphi$,  and $P_x(k+1,k+2,\psi)=P_x(k+1,k+2,\varphi)$. 
However, we have $ \beta\in (\psbar(v_{i-1})\cap \psbar(x)) \cap [1,k] \ne \emptyset$, contradicting the choice of $\varphi$. 
Thus $c_\varphi=0$, completing the proof. 
\qed 

Let $G$ be a graph and 
$k=\min\{\delta(G)-1, \lfloor \rho_c(G) \rfloor \}$. A subset $U$ of $V(G)$
is odd if $|U| \ge 3$ is odd. 
An odd  set $U$ of $G$
is \emph{optimal} (with respect to $k$) if $e^+(U)= k (|U|+1)/2$.   
 For an optimal set $U$ of $G$,    since $2e^+(U)=\sum_{v\in U}d(v)+e(U,V(G)\setminus U)$,
we get $k (|U|+1)=\sum_{v\in U}d(v)+e(U,V(G)\setminus U) \ge (k+1)|U|+e(U,V(G)\setminus U)$ with equality holds if $\sum_{v\in U}d(v)= (k+1)|U|$. Thus 
\begin{eqnarray}
	k &\ge& |U|+e(U,V(G)\setminus U) \quad \text{and}  \nonumber \\ 
	&& \label{eqn1} \\
	k &= &|U|+e(U,V(G)\setminus U) \quad \text{if $\sum_{v\in U}d(v)= (k+1)|U|$.} \nonumber 
\end{eqnarray} 
We have the following property for optimal sets of $G$. 
\begin{LEM}\label{lemma:disjoint-optimal sets}
	Let $G$ be a graph with  $k=\min\{\delta(G)-1, \lfloor \rho_c(G) \rfloor \}$ and $x\in V(G)$. Suppose   $U$ is  an optimal set of $G$ such that $x\in U$ and $|U|$
	is minimum among the sizes of all optimal sets containing $x$. 
	Then for any  optimal set $U'$ of $G$ with $U\not\subseteq U'$, we have $x\not\in U\cap U'$.  
	\end{LEM}
\pf  Suppose to the contrary that $x\in U\cap U'$. 
Let 
$$
L=U\setminus U', \quad M=U\cap U',  \quad  R=U'\setminus U, \quad \text{and} \quad W=V(G)\setminus (U\cup U'). 
$$
Since $U\not\subseteq U'$ and $U\cap U' \ne \emptyset$, we have $L, M  \ne \emptyset$.  As $|U|$ is  minimum among the sizes of all optimal sets containing $x$, 
 $x\in U'$ and $U'\ne U$, we have $U'\not\subseteq U$. Thus $R\ne \emptyset$ as well. 
By counting the edges within distinct parts, we have 
\begin{eqnarray*}
	e^+(U\cup U')&=&e(L)+e(M)+e(R)+e(L,M)+e(M,R)+e(L,R)+e(L,W)+ \\
	&&e(M,W)+e(R,W), \\
	e^+(U) &=&e(L)+e(M)+ e(L,M)+e(M,R)+e(L,R)+e(L,W)+e(M,W),\\
	e^+(U') &=&e(R)+e(M)+e(L,M)+e(M,R)+e(L,R)+e(R,W)+e(M,W),\\
	e^+(M) &=& e(M)+e(L,M)+e(M,R)+e(M,W), \\
	e^+(L) &=& e(L)+e(L,M)+e(L,R)+e(L,W) ,\\
	e^+(R) &=& e(R)+e(M,R)+e(L,R)+e(R,W). 
\end{eqnarray*}
Therefore,  $$e^+(U\cup U')=e^+(U)+e^+(U')-e^+(M)-e(L,R).$$
When $|M|$ is odd, if $|M|=1$, then $e^+(M) \ge \delta(G) \ge k+1=k (|M|+1)/2+1$. 
Otherwise $|M|\ge 3$.   Since  $|M|<|U|$ and $x\in M\subseteq U$, we know that $M$
is not optimal by the choice of $U$. Thus $e^+(M)\ge k (|M|+1)/2+1$.

Suppose first that $|M|$ is odd and so $|U\cup U'|$ is odd. Then 
\begin{eqnarray*}
	e^+(U\cup U')&=&e^+(U)+e^+(U')-e^+(M)-e^+(L,R) \\
	& \le & k (|U|+1)/2 +k (|U'|+1)/2 -(k (|M|+1)/2+1)-e^+(L,R) \\
	&= & k (|U\cup U'|+1)/2-1-e^+(L,R) < k (|U\cup U'|+1)/2,  
\end{eqnarray*}
a contradiction to the assumption that $ \lfloor \rho_c(G) \rfloor  \ge k$. 

Thus we assume that $|M|$ is even. Then $|L|$ and $|R|$ are odd. 
Again we have $e^+(L)\ge k (|L|+1)/2$ and $e^+(R)\ge k (|R|+1)/2$ by the assumption that $k=\min\{\delta(G)-1, \lfloor \rho_c(G) \rfloor \}$. 
As $2e(M)+e(L,M)+e(M,R)+e(M,W)=\sum_{x\in M}d(x) \ge (k+1)|M|$, we get 
\begin{eqnarray*}
	e^+(U)+e^+(U') &=& e^+(L)+e^+(R)+2e(M)+e(L,M)+e(M,R)+2e(M,W) \\
	& \ge &k (|L|+1)/2+k (|R|+1)/2+(k+1)|M|+e(M,W) \\
	&\ge & k(|L|+1)/2+k (|R|+1)/2 +k|M|/2+k|M|/2+|M| \\
	&=& k (|U|+1)/2+k (|U'|+1)/2+|M| \ge k(|U|+1)/2+k (|U'|+1)/2+1, 
\end{eqnarray*}
a contradiction to  the assumption that both $U$ and $U'$ are optimal. 
\qed

\section{Proof of Theorem~\ref{thm:main}}

Let $V=V(G)$ and $E=E(G)$, and  $k=\min\{\delta(G)-1, \lfloor \rho_c(G) \rfloor \}$. Then  $\delta(G)\ge k+1$
and for any odd $U\subseteq V(G)$, we have $e^+_{G}(U) \ge k (|U|+1)/2$. 
Recall that an odd $U\subseteq V(G)$
is  optimal if $e^+_{G}(U)= k (|U|+1)/2$. 
For a  vertex $x\in V$ with $d_G(x) \ge k+2$, we apply the following operation:
\begin{itemize}
	\item If $x$ is not contained in any optimal set of $G$, then we \emph{split off} an edge $xy$ that is incident with $x$ from $x$ (i.e.,  add a new vertex $x'$, delete $xy$ but add the edge $yx'$); 
	\item If $x$ is contained in an optimal set of $G$, we let $U$ be an optimal set containing $x$ with minimum size.  
	Let $y\in U$ with $xy\in E$, and then we split $xy$ off from $x$.  Note that the vertex $y$ exists 
	as $e_G(x,V\setminus U) \le k-|U|$ by Equation~\eqref{eqn1}. 
\end{itemize}

Let $H$ be the resulting graph.  We claim that for any odd set  $U \subseteq V$ of $H$, 
we still have  $e^+_H(U) \ge k (|U|+1)/2$. 
Suppose to the contrary that there exists $U'\subseteq V$ such that $e^+_H(U') \le k (|U'|+1)/2-1$. As 
we only split off one edge from $x$ in $G$ to get $H$,  it follows that $x\in U'$
and $y\not\in U'$,  
$e^+_H(U')=k (|U'|+1)/2-1$ and $e^+_G(U')=k (|U'|+1)/2$. Thus $U'$ is  optimal in $G$. 

It must be the case that $x$ is contained in an optimal set of $G$ when we obtained $H$ from $G$. 
We let $U$
be the  optimal set of $G$ for which $x,y\in U$ and we split  $xy$ off from $x$ when we did the splitting off operation. 
Note that $U$ was chosen to be an optimal set containing $x$ with minimum size. 
 Now we have $x\in U\cap U'$, $y\in U\setminus U'$ and so $U\not\subseteq U'$.  However, this shows a contradiction to Lemma~\ref{lemma:disjoint-optimal sets}. 

Thus we can repeat the operation on $x$ for all vertices of $V$ that has degree at least $k+2$
in $H$. Still denote the resulting graph by $H$ and we can assume now  $d_H(v)=k+1$ for any $v\in V$.  
By the operation, any $v\in V(H)\setminus V$ satisfies $d_H(v)=1$.  
As $E(H[V]) \subseteq E$ and every edge from $e_H(V, V(H)\setminus V)$ 
corresponds to an edge of $E$, it suffices to show that $H$
has $k$ disjoint edge sets that each saturate  $V$.

For any odd  $U\subseteq V$ of $H$, we have $e_H(U)+e^+_H(U)=2e_H(U)+e_H(U,V(H)\setminus U)=(k+1)|U|$.
Thus 
\begin{numcases}{e_H(U) \le }
	(k+1)|U|-k(|U|+1)/2-1  & \text{if $U$ is not optimal}; \nonumber  \\ 
	= k(|U|-1)/2+|U|-1  =(k+2)(|U|-1)/2& \nonumber \\ 
	&\label{eqn2}\\
	k(|U|-1)/2+|U|= (k+2)(|U|-1)/2 +1 & \text{if $U$ is  optimal}.   \nonumber  
\end{numcases}

By Lemma~\ref{lemma:disjoint-optimal sets}, all minimum optimal 
sets contained in $V$ are vertex-disjoint (by the lemma, any vertex of a  minimum optimal set is contained only in one minimum optimal set). 
If exist,  
let  $U_1, U_2, \ldots, U_t$  be all the minimum 
 optimal  sets of $H$  that are contained in $V$, where  $t\ge 1$ is an integer.

For each $i\in [1,t]$, we delete an edge $x_iy_i$ from $H[U_i]$. Denote the resulting graph 
by $H_1$.  We claim that $\chi'(H_1[V]) =k+2$. 
For otherwise, $\chi'(H_1[V]) =s+1 \ge k+3=\Delta(H)+2$ for some integer $s$.   
As $H_1[V]$ contains an edge-chromatic critical subgraph (a connected graph  whose  every edge  is critical) with chromatic index being $s+1$, 
Lemma~\ref{G-S-C} implies that  there is an edge $e\in E(H_1[V])$ and a subgraph $J\subseteq H_1[V]$
such that $J-e$ is $s$-dense or $e_{H_1}(J-e)=s(|V(J)|-1)/2$. 
This implies that $e_{H_1}(V(J))=(k+2)(|V(J)|-1)/2+1$ by~\eqref{eqn2} and so $s=k+2$. 
As  $|V(J)| \ge 3$
and is odd, we know that $V(J)$ is an optimal  set in $H_1$ by~\eqref{eqn2}. 
We let $U^* \subseteq V$ be a minimum optimal  set of $H_1$. 
Thus $e_{H_1}(U^*) = (k+2)(|U^*|-1)/2+1$. 
Since $e_{H}(U^*) \le (k+2)(|U^*|-1)/2+1$ by~\eqref{eqn2}, we know that 
 $U^*$ is an optimal set in $H$ as well.   Thus $U^*=U_i$
 for some $i\in [1,t]$ by Lemma~\ref{lemma:disjoint-optimal sets}. As an edge $x_iy_i$ was removed from $H[U_i]$
 in getting $H_1$, we have $e_{H_1}(U^*)= (k+2)(|U^*|-1)/2$, a contradiction 
 the previous assumption that  $e_{H_1}(U^*) = (k+2)(|U^*|-1)/2+1$. Thus $\chi'(H_1[V]) =k+2$.  
 As vertices in $V(H_1)\setminus V$ have degree 1 in $H_1$, we have $\chi'(H_1)=k+2$. 
 
For each $i\in [1,t]$, as $e_{H_1}(U_i)=e_{H}(U_i)-1=(k+2)(|U_i|-1)/2$, we know 
that for any $\varphi\in \CC^{k+2}(H_1)$, the colors on 
the edges in $E_{H_1}(U_i, V(H_1)\setminus U_i)$ under $\varphi$ are all distinct. 
Thus the graph $H_2$ obtained from $H_1$ by contracting each $U_i$ into a single vertex $u_i$
for each $i\in [1,t]$ is edge $(k+2)$-colorable.

We claim that $d_{H_2}(u_i) \le k/2$ for each $i\in [1,t]$. 
Suppose, without loss of generality, that $|U_1| \le |U_2| \le \ldots \le |U_t|$. 
Then by~\eqref{eqn1}, we have $e_{H}(U_1, V(H)\setminus U_1) \ge e_{H}(U_2, V(H)\setminus U_2) \ge \ldots  \ge e_{H}(U_t, V(H)\setminus U_t)$. 
Since $H_1$ was obtained from $H$ by deleting one edge within each $U_i$, we have  $d_{H_2}(u_i)=e_{H}(U_i, V(H)\setminus U_i)$. Thus 
$d_{H_2}(u_1) \ge d_{H_2}(u_2) \ge \ldots \ge d_{H_2}(u_t)$. 
It then suffices to show that $d_{H_2}(u_1) \le k/2$, or equivalently $e_H(U_1, V(H)\setminus U_1) \le k/2$.  
As $(k+2)(|U_1|-1)+2=2e_H(U_1) $, when the maximum multiplicity of 
$G$ is at most 2, we have $2e_{H}(U_1) \le 2(|U_1|-1)|U_1|$ and so $k+2 \le 2|U_1|$. 
Now by~\eqref{eqn1} that $k=|U_1|+e_H(U_1, V(H)\setminus U_1)$, we get $e_H(U_1, V(H)\setminus U_1) =k- |U_1| \le k-k/2=k/2$.
When $k\le 6$, then as $|U_1| \ge 3$, $k=|U_1|+e_H(U_1, V(H)\setminus U_1)$ from~\eqref{eqn1} implies that $e_H(U_1, V(H)\setminus U_1) \le k/2$.  
Therefore $d_{H_2}(u_1) \le k/2$ and thus $d_{H_2}(u_i) \le k/2$ for each $i\in [1,t]$.

By Lemma~\ref{lemma:special-coloring}, $H_2$ has an edge $(k+2)$-coloring $\varphi$ 
satisfying the following two properties: 
(1) the color $k+2$ is missing at every vertex in $\{u_1, \ldots, u_t\}$;  and (2)
if $k+1\in \varphi(u_i)$ for some $i\in [1,t]$, then $P_{u_i}(k+1,k+2, \varphi)$ does not end 
	at any vertex from $\{u_1, \ldots, u_t\}\setminus \{u_i\}$.

We extend the coloring $\varphi$ of $H_2$ to a coloring $\psi$ of $H_1$, for each $i\in [1,t]$, 
by giving an edge  $(k+2)$-coloring  $\varphi_i$ of each $H_1[U_i]$ and then permute the colors  of $\varphi_i$ so 
that if $wz\in E_{H_1}(U_i, V(H_1)\setminus U_i)$ with $w\in U_i$, then  the color $\varphi(u_i z)$
is missing at $w$ under $\varphi_i$. This is possible by Lemma~\ref{k-dense subgraph}.  
In particular, since the class classes of $\varphi_i$ are  near perfect matchings 
such that each vertex $w\in U_i\setminus \{x_i,y_i\}$ is missed by  exactly $1+e_{H_1}(w, V(H_1)\setminus U_i)$ of them 
and both $x_i$ and $y_i$ are missed by $2+e_{H_1}(w, V(H_1)\setminus U_i)$ of them, 
we can assume  that a  near perfect matching  of $H_1[U_i]$ that misses the vertex $x_i$ is colored by the color $k+2$ under $\psi$. 
Thus we have an edge $(k+2)$-coloring $\psi$ of $H_1$ such that 
\begin{enumerate}
	\item  the color $k+2$ is missing on edges in $E_{H_1}(U_i, V(H_1)\setminus U_i)$ for each $i\in [1,t]$; 
	\item If  the color $k+1$  presents on an edge from $ E_{H_1}(U_i, V(H_1)\setminus U_i)$ for some $i\in [1,t]$,  
	then the corresponding $(k+1,k+2)$-chain 
	including that edge  ends at a vertex in $V(H_1)\setminus (\bigcup_{ 1\le i\le t} U_i)$. In this case, all vertices in $U_i$
	present the color $k+1$ under $\psi$; 
	\item The $(k+2)$ color class   misses $x_i$ for each $i\in [1,t]$.  
\end{enumerate} 

We let $M_1, \ldots, M_{k+2}$ be the color classes of $\psi$  corresponding to the colors $1, \ldots, k+2$ respectively.  
We will add edges from $M_{k+1}\cup M_{k+2}$ to each of $M_1, \ldots, M_k$ 
if necessary to modify them into $k$ edge-disjoint edge sets  of $H_1$ that each saturate $V$.  
To do so, consider $D'$ the subgraph  of $H_1$ induced on  $M_{k+1}\cup M_{k+2}$. 
As $\Delta(D_1) \le 2$, the graph can be decomposed into vertex-disjoint cycles and paths. 
We orient $D_1$ such that the cycles are directed cycles and 
 the paths are directed paths. Moreover, 
if a path has $x_i$ as one endvertex for some $i\in [1,t]$, then 
the path is oriented towards $x_i$. Note  that if a path has $x_i$ has one of 
its endvertex, then its another endvertex is a distinct vertex 
 not in $\{u_1, \ldots, u_t\}$ by our choice of the coloring $\psi$. 
Denote the orientation of $D'$ by $D$. 

Each vertex from $V\setminus\{x_i,y_i : i\in [1,t]\}$ has degree $k$ in $H_1$, and so it  is missed by at most one 
of $M_1, \ldots, M_k$. If it is missed by exactly one of $M_1, \ldots, M_k$, then it has degree 2 in $D'$. 
The vertex $x_i$ for each $i\in [1,t]$ has degree at most 1 in $D'$ by our choice of $\psi$. 
Thus each vertex $w\in V\setminus\{y_i : i\in [1,t]\}$  is missed by at most one matching 
from $M_1, \ldots, M_k$. If $w$ is missed by exactly one of $M_1, \ldots, M_k$, by our construction, 
there is in $D$ an arc $zw$.  We then add $zw$ to the matching that misses $w$. 
Now each vertex $y_i$ for $i\in [1,t]$  is missed by at most two matchings 
from $M_1, \ldots, M_k$. If $y_i$ is missed by exactly one matching from 
$M_1, \ldots, M_k$, then we add the edge $x_iy_i$ to the matching that misses $y_i$. 
 If $y_i$ is missed by exactly two matchings from $M_1, \ldots, M_k$, 
then it has degree 2 in $D'$ and so there is in $D$ an arc $zy_i$. Then we add $zy_i$ to one of the matching that misses $y_i$
and add $x_iy_i$ to the other matching that misses $y_i$.  After those modifications, we have obtained $k$
edge-disjoint edge sets of $H$ based on $M_1, \ldots, M_k$  and each of those sets saturates  $V$. The proof is then complete. 
\qed

\section*{Acknowledgment}
Guantao Chen was supported by NSF grant DMS-2154331 and Songling Shan was supported by NSF grant DMS-2153938.

\bibliographystyle{plain}
\bibliography{SSL-BIB_08-19}

\begin{thebibliography}{1}

\bibitem{1906.06458}
Y.~Cao, G.~Chen, G.~Ding, G.~Jing, and W.~Zang.
\newblock On gupta's co-density conjecture.
\newblock {\em \arxiv{1906.06458}}, 2019.

\bibitem{1901.10316}
G.~Chen, G.~Jing, and W.~Zang.
\newblock Proof of the goldberg-seymour conjecture on edge-colorings of
  multigraphs.
\newblock {\em \arxiv{1901.10316}}, 2019.

\bibitem{MR354429}
M.~K. Goldberg.
\newblock Multigraphs with a chromatic index that is nearly maximal.
\newblock {\em Diskret. Analiz}, (23):3--7, 72, 1973.

\bibitem{Gupta}
R.~P. Gupta.
\newblock On decompositions of a multigraph into spanning subgraphs.
\newblock {\em Bull. Amer. Math. Soc.}, (80):500--502, 1974.

\bibitem{MR0491300}
R.~P. Gupta.
\newblock On the chromatic index and the cover index of a multigraph.
\newblock In {\em Theory and applications of graphs ({P}roc. {I}nternat.
  {C}onf., {W}estern {M}ich. {U}niv., {K}alamazoo, {M}ich., 1976)}, Lecture
  Notes in Math., Vol. 642, pages 204--215. Springer, Berlin, 1978.

\bibitem{MR532981}
P.~D. Seymour.
\newblock On multicolourings of cubic graphs, and conjectures of {F}ulkerson
  and {T}utte.
\newblock {\em Proc. London Math. Soc. (3)}, 38(3):423--460, 1979.

\end{thebibliography}

\end{document}